\documentclass[12pt]{article}
\usepackage{amssymb,amsmath}
\def \R{{\hbox{\vrule width 0.6pt height 6.8pt depth -.2pt\kern-0.2pt
R}}}
\def \P{{\hbox{\vrule width 0.6pt height 6.8pt depth -.2pt\kern-0.2pt
P}}}

\def \no { \noindent}
\def \er {\mathbb R}

\def\grad{\nabla}
\def\div{\, \mbox{div}\,  }

\def\MOD#1{{|\kern -.16em |\kern -.16em | #1 | \kern -.16em |\kern
 -.16em |}}
\def \epsilon {\varepsilon}

\def\ds{\displaystyle}
\newtheorem{theo}{\bf THEOREM}[section]
\newtheorem{lem}[theo]{\bf LEMMA}
\newtheorem{pro}[theo]{\bf PROPOSITION}
\newtheorem{cor}[theo]{\bf COROLLARY}

\numberwithin{equation}{section}

\def\Box{\hfill\rule{2.5mm}{2.5mm}}
\begin{document}
\bibliographystyle{siam}

\title{
A Lyapunov functional and  blow-up results for a class of perturbed  semilinear wave
equations}
\author{M. Hamza\\
{\it \small Facult\'e des Sciences de Tunis}\\
H. Zaag\footnote{This author is supported by a grant from the french Agence Nationale de la Recherche, project ONDENONLIN, reference ANR-06-BLAN-0185.} \\
{\it \small CNRS LAGA Universit\'e Paris 13}}

%\author{ M. Hamza and H. Zaag}

\maketitle

\begin{abstract}
\no We  consider in this paper some class of perturbation for the semilinear wave equation with subcritical (in the conformal transform sense) power nonlinearity.
We first derive a
Lyapunov functional in similarity
 variables and then use it to derive the blow-up rate. Though the result is similar to the unperturbed case in its statements, this is not the case of our method, which is new up to our knowledge.

\end{abstract}


\vspace{0.5cm}


 \vspace{0.4cm}

\maketitle \maketitle
\section{\bf Introduction}
This paper is devoted to the study of blow-up solutions for the following semilinear  wave
equation:
\begin{equation}
\label{1} \left\lbrace
\begin{array}{l}
 u_{tt}=
\Delta u+ |u|^{p-1} u+f(u)+g(u_t),\qquad  (x,t)\in \er^N\times \er_+^* \\
\\
(u(x,0),u_t(x,0))=(u_0(x),u_1(x))\in H^{1}_{loc,u}(\er^N)\times L^{2}_{loc,u}(\er^N),\\
\end{array}
\right.
\end{equation}
where $p>1$,
$f$ and $g$  are locally Lipschitz-continuous  satisfying
the following condition
$$(H_f)\qquad |{f(x)}|\le M(1+|x|^q) \qquad{{\textrm {with}}}\ \ (q<p,\ \ M>0),$$
$$(H_g)\qquad \qquad \qquad |{g(x)}|\le M(1+|x|),\qquad\qquad$$
 and $L^{2}_{loc,u}(\er^N)$ and $H^{1}_{loc,u}(\er^N)$ are the spaces defined by
\begin{equation*}
L^{2}_{loc,u}(\er^N)=\{u:\er^N\rightarrow \er/ \sup_{a\in \er^N}(\int_{|x-a|\le 1}|u(x)|^2dx)<+\infty \},
\end{equation*}
and
\begin{equation*}
H^{1}_{loc,u}(\er^N)=\{u\in L^{2}_{loc,u}(\er^N),|\grad u|\in L^{2}_{loc,u}(\er^N) \}.
\end{equation*}
We assume in addition  that $$1<p<p_c\equiv 1+\frac{4}{N-1}.$$
\no
The  Cauchy problem  of  equation ({\ref{1}}) is  wellposed
in $H^{1}_{loc,u}\times L^{2}_{loc,u}$. This follows from the finite
speed of propagation and the the wellposdness in $H^{1} \times L^{2}$, valid whenever
$ 1< p<1+\frac{4}{N-2}$ .
The existence of blow-up   solutions for  the  associated  ordinary differential equation  of  (\ref{1})  is a classical result. By using  the
finite
speed of propagation, we conclude that there exists a  blow-up solution $u(t)$ of (\ref{1}).
In this paper, we consider  a blow-up solution $u(t)$ of (\ref{1}), we define (see for example Alinhac \cite{A} and \cite{A1})  $\Gamma$ as the graph of a function $x\rightarrow T(x)$ such that $u$ cannot be extended beyond the set
$$D_u=\{(x,t)\ \ \big |t<T(x)\}.$$
The set $D_u$ is called the maximal influence domain of $u$. Moreover, from the finite speed of propagation, $T$ is a $1$-Lipschitz function.
Let $\overline{T}$ be the minimum of $T(x)$ for all $x\in \er^N$. The time
$\overline{T}$  and the graph $\Gamma$  are called (respectively) the blow-up
time and the blow-up graph of $u$. \\
Let us first  introduce the following non degeneracy condition for $\Gamma$. If we introduce for all $x\in\er^N$, $t\le T(x)$ and $\delta>0$, the cone

\begin{equation} \label{cone}
C_{x,t,\delta}=\{(\xi,\tau)\neq(x,t)|0\le \tau\le t-\delta|\xi-x|\},
\end{equation}
then our non degeneracy condition is the following: $x_0$ is a non characteristic point if
\begin{equation} \label{cone1}
\exists\delta_0=\delta_0(x_0)\in (0,1) \ \textrm{ such that}\  u \ \textrm{ is defined on }\ C_{x_0,T(x_0),\delta_0}.
\end{equation}
%We have in addition
%$$\|u(t)|_{H^{1}_{{loc,u}}}+\|\partial_tu(t)|_{L^{2}_{{loc,u}}}\rightarrow +\infty \ \textrm{as} \ \ t\rightarrow \overline{T}.$$
We aim at studying  the growth estimate  of $u(t)$ near the space-time blow-up graph.

\no In the case $(f, g)\equiv (0,0)$, equation (\ref{1}) reduces to the semilinear wave equation:
\begin{equation}\label{mu=0}
 u_{tt}=\Delta u+ |u|^{p-1} u,\qquad  (x,t)\in \er^N\times [0,\overline{T}).
\end{equation}

\no  Merle and Zaag   in \cite{MZ3} (see also    \cite{MZ1} and \cite{MZ2}) have proved, that if $u$ is a solution of (\ref{mu=0})
with blow up graph
$\Gamma:\{x\rightarrow T(x)\}$, then for all $x_0 \in \er^N$  and $t\in [\frac34T(x_0),T(x_0)]$,
  the growth estimate near the space-time blow-up graph  satisfies
\begin{eqnarray*}
 (T(x_0)-t)^{\frac{2}{p-1}}\frac{\|u(t)\|_{L^2(B(x_0,\frac{T(x_0)-t}{2}))}}{ (T(x_0)-t)^{\frac{N}{2}}}\qquad \qquad\qquad\qquad\qquad\qquad \qquad\qquad\qquad\\%\qquad\\
 + (T(x_0)-t)^{\frac{2}{p-1}+1}\!\Big (\frac{\|u_t(t)\|_{L^2(B(x_0,\frac{T(x_0)-t}{2}))}}{ (T(x_0)-t)^{\frac{N}{2}}}+
 \frac{\|\grad u(t)\|_{L^2(B(x_0,\frac{T(x_0)-t}{2}))}}{ (T(x_0)-t)^{\frac{N}{2}}}\Big )\le K,
\end{eqnarray*}
where the constant $K$  depends only on $N,p,$ and on  an upper bound on $T(x_0)$, $\frac1{T(x_0)}$
and the initial data in $ H^{1}_{loc,u}(\er^N)\times L^{2}_{loc,u}(\er^N)$.
If in addition $x_0$ is non characteristic (in the sense (\ref{cone1})), then for all $t\in [ \frac{3T(x_0)}4, T(x_0)]$,

\begin{eqnarray*}
0<\varepsilon_0(N,p)\le (T(x_0)-t)^{\frac{2}{p-1}}\frac{\|u(t)\|_{L^2(B(x_0,{T(x_0)-t}))}}{ (T(x_0)-t)^{\frac{N}{2}}}
\qquad \qquad\qquad\qquad\qquad
\nonumber\\
 + (T(x_0)-t)^{\frac{2}{p-1}+1}\Big (\frac{\|u_t(t)\|_{L^2(B(x_0,{T(x_0)-t}))}}{ (T(x_0)-t)^{\frac{N}{2}}}+
 \frac{\|\grad u(t)\|_{L^2(B(x_0,{T(x_0)-t}))}}{ (T(x_0)-t)^{\frac{N}{2}}}\Big )\le K,
\end{eqnarray*}
where the constant $K$  depends only on $N,p,$ and on  an upper bound on $T(x_0)$, $\frac1{T(x_0)}$, $\delta_0(x_0)$
and the initial data in $ H^{1}_{loc,u}(\er^N)\times L^{2}_{loc,u}(\er^N)$.

\no Following this blow-up rate estimate, Merle and Zaag addressed the question of the asymptotic
behavior of $u(x,t)$ near $\Gamma$
in one space dimension.

\no More precisely, they proved in \cite{MZ4}
 and \cite{MZ5} that the set of non charecteristic points ${\cal R}\subset \er$ is open
and that $x\mapsto T(x)$ is of class $C^1$ on ${\cal R}$. They also described the blow-up profile of $u$  near $(x_0,T(x_0))$ when $x_0\in {\cal R}$.

\no In \cite{MZarxiv}, they proved that $S=\er\backslash {\cal R}$ has an empty interior and that $\Gamma$ is a corner  of angle $\frac{\pi}{2}$ near any $x_0\in S$. They also showed that
$u(x,t)$ decomposes in a sum  of decoupled solitons near $(x_0,T(x_0))$.

\bigskip
\no Our aim in this work is to generalize the blow-up rate estimate obtained for   equation (\ref{mu=0})
 in \cite{MZ1} and \cite{MZ3} in the subcritical case $(p<p_c)$
to  equation (\ref{1}).

\no One may think that such a generalization   is straightforward and only technical.
 In fact, that opinion may be valid for all the steps, except for the very first one, that is, the existence of a  Lyapunov functional
in  similarity variables which  is far from being trivial. That functional is our main contribution.
The existence of the Lyapunov functional is a crucial step towards the derivation of blow-up results for equation (\ref{1})

\no As in \cite {MZ1} and \cite{MZ2}, we want to write
the solution $v$ of the associate ordinary differential equation of (\ref{1}).
It is clear that $v$ is given by
\begin{equation}
v''= v^{p} +f(v)+g(v'),\qquad v({T})=+\infty,
\end{equation}
and satisfies:
\begin{equation}\label{kappa}
v(t)\sim \frac{\kappa}{ ({T}-t)^{\frac{2}{p-1}}}\ \  \textrm{as}\ t\rightarrow {T},
\ \  \textrm{where}\ \kappa=\Big (\frac{2p+2}{(p-1)^2}\Big )^{\frac{1}{p-1}}.
\end{equation}
For this reason, we define
for all $x_0\in \er^N$, $0< T_0\le T_0(x_0)$,
 the following similary transformation introduced
in Antonini and Merle \cite{AM} and used in \cite{MZ1},\cite{MZ2} and \cite{MZ3}:
\begin{equation}\label{scaling}
y=\frac{x-x_0}{T_0-t},\qquad s=-\log (T_0-t),\qquad u(x,t)=\frac{1}{(T_0-t)
^{\frac{2}{p-1}}}w_{x_0,T_0}(y,s).
\end{equation}
The function $w_{x_0,T_0}$  (we write $w$ for simplicity)
satisfies the following equation for all $y\in B$ and $s\ge -\log T_0$:
\begin{eqnarray}\label{B}
 w_{ss}&=&\frac{1}{\rho}\div(\rho \grad w-\rho(y.\grad w)y)
-\frac{2p+2}{(p-1)^2}w+|w|^{p-1}w-\frac{p+3}{p-1}w_s-2y.\grad w_{s}\nonumber\\
&&+e^{-\frac{2ps}{p-1}}f\Big(e^{\frac{2s}{p-1}}w\Big) +
e^{-\frac{2ps}{p-1}}g\Big(e^{\frac{(p+1)s}{p-1}}(w_s+y.\grad w+\frac{2}{p-1}w)\Big),\qquad
\end{eqnarray}
where
$\rho=(1-|y|^2)^{\alpha}$ and
$\alpha={\frac{2}{p-1}-\frac{N-1}{2}}>0.$

\no In the new set of variables $(y,s),$
the behavior of $u$ as $t \rightarrow T_0$
is equivalent to the
behavior of $w$ as $s \rightarrow +\infty$.

\no {\bf{Remark:}}\\
We remark that the corresponding  terms of the fonctions $f(u)$ and $g(u)$  in the problem (\ref{B}) satisfy the following inequalities, for all $s\ge 0$,
\begin{eqnarray*}
e^{-\frac{2ps}{p-1}}\left|f\left(e^{\frac{2s}{p-1}}w\right)\right| &\le&
CMe^{-\frac{2ps}{p-1}}
+CMe^{-\frac{2(p-q)s}{p-1}}|w|^q\\
&\le& CMe^{-\frac{2(p-q)s}{p-1}}
+CMe^{-\frac{2(p-q)s}{p-1}}|w|^p
,
\end{eqnarray*}
and
\begin{eqnarray*}
e^{-\frac{2ps}{p-1}}g\left|\left(e^{\frac{(p+1)s}{p-1}}(w_s+y.\grad w+\frac{2}{p-1}w)\right)\right|\quad\quad\quad\quad\quad\quad\quad\quad\quad\quad\quad\quad\quad\quad\quad\quad\\
\le
CMe^{-\frac{2ps}{p-1}}
+CMe^{-s}\left|w_s+y.\grad w+\frac{2}{p-1}w\right|.
\end{eqnarray*}
%\begin{eqnarray}
%e^{-\frac{2ps}{p-1}}\Big|f\Big(e^{\frac{2s}{p-1}}w\Big)\Big| &\le&
%CMe^{-\frac{2ps}{p-1}}
%+CMe^{-\frac{2(p-q)s}{p-1}}|w|^q\le CMe^{-\frac{2(p-q)s}{p-1}}
%+CMe^{-\frac{2(p-q)s}{p-1}}|w|^p
%,\nonumber
%\end{eqnarray}
%\begin{eqnarray}
%e^{-\frac{2ps}{p-1}}g\Big|\Big(e^{\frac{(p+1)s}{p-1}}(w_s+y.\grad w+\frac{2}{p-1}w)\Big)\Big|&\le&
%CMe^{-\frac{2ps}{p-1}}
%+CMe^{-s}\Big|w_s+y.\grad w+\frac{2}{p-1}w\Big|.\nonumber
%\end{eqnarray}
For this reason, we can see that in the variable $(y,s)$ the problem (\ref{B}) is a perturbation of the particular case  where $(f,g)\equiv(0,0)$, when $s\rightarrow +\infty$.

\vspace{0.3cm}
\no
The equation (\ref{B}) will be studied in the space $\cal H$
$${\cal H}=\Big \{(w_1,w_2), |
\displaystyle\int_{B}\Big ( w_2^2
+|\grad w_1|^2(1-|y|^2)+w_1^2\Big )\rho {\mathrm{d}}y<+\infty \Big \}.$$
In the whole paper,   we  denote
$\ds{F(u)=\int_0^uf(v)dv}.$\\
In the case $(f,g)\equiv(0,0)$,  Antonini and Merle \cite{AM} proved that
\begin{eqnarray}
\label{f}
E_0(w)=\displaystyle\int_{B}\Big ( \frac{1}{2}w_s^2
+\frac{1}{2}|\grad w|^2-\frac{1}{2}(y.\grad w)^2+\frac{p+1}{(p-1)^2}w^2
-\frac{1}{p+1}|w|^{p+1}\Big )\rho {\mathrm{d}}y,
\end{eqnarray}
is a Lyapunov functional for equation (\ref{B}).
 When $(f,g) \not\equiv (0,0)$, we introduce
\begin{eqnarray}
\label{f1}
 H(w)&=&E(w)e^{\frac{ p+3}{2\gamma} e^{-\gamma s}}+\theta e^{-2\gamma s},
\end{eqnarray}
where $\theta$  is a sufficiently large constant that will be determined later,
\begin{eqnarray}
\label{f2}
 E(w)&=&E_0(w)+I(w)+J(w),\quad
  I(w)=- e^{-\frac{2(p+1)s}{p-1}}\displaystyle\int_{B}F(e^{\frac{2}{p-1}s}w)\rho {\mathrm{d}}y,\nonumber\\
  {\textrm {and}} \quad
 J(w)&=& -e^{-\gamma s}\displaystyle\int_{B}ww_s\rho {\mathrm{d}}y \qquad {\textrm {with}} \qquad \gamma=\min(\frac12,
 \frac{p-q}{p-1}
 )>0.%\qquad\qquad\qquad\qquad
\end{eqnarray}
\no We now claim that the functional $H(w)$ is a decreasing function of time for equation
 $(\ref{B})$, provided that $s$ is   large enough.

\bigskip
\no Here we announce our main result.

\begin{theo}(Existence  of a Lyapunov functional for equation ({\ref{B}}))\\
\label{lyap}
Let $N,p,q$  and $M$ be fixed. There exists  $S_0=S_0(N,p,q,M)\in \er$
such that, for all $ s_0\in \er$ and $w$ solution of
equation (\ref{B}) satisfying  $(w,w_s)\in
 {\cal
  C}([s_0,+\infty),{\cal H})$, it holds that $H$ satisfies the following inequality, for all $s_2>s_1\ge \max(s_0,S_0)$,
\begin{eqnarray}\label{L1}
H(w(s_2))-H(w(s_1))
&\le&  -\alpha \int_{s_1}^{s_2}\int_{B}w_s^2(y,s)\frac{\rho}{1-|y|^2}dy ds, \qquad \qquad
\end{eqnarray}
where $\alpha=\frac{2}{p-1}-\frac{N-1}2$.
\end{theo}

\no {\bf{Remarks:}}
\begin{enumerate}
\item One may wander why  we take only sublinear perturbations in $u_t$ (see hypothesis ($H_g$)).
It happens that  any superlinear  terms in $u_t$ generates in similarity variables $L^r$ norms of $w_s$ and $\grad w$, where $r>2$, hence,
non controllable by the terms in the Lyapunov functional $E_0(w)$ (\ref{f}) of the non perturbed equation (\ref{mu=0}).
\item Our method  breaks   down in the critical
 case $p=p_c$, since in the energy estimates in similarity variables, the perturbations terms are integrated on the whole unit ball, hence, difficult to control with the dissipation of the non perturbed equation (\ref{mu=0}), which degenerates to the boundary of the unit ball.
\end{enumerate}
As we said earlier, the existence of this Lyapunov functional (and a blow-up criterion for equation (\ref{B}) based in $H$, see Lemma 2.3 below)  are a crucial step in the derivation of the blow-up rate for equation (\ref{1}).
Indeed, with the  functional $H$ and some more work, we are able to adapt the analysis performed in \cite{MZ3} for equation (\ref{mu=0}) and get the following result:
\begin{theo}(blow-up rate  for equation ({\ref{1}}))\\
\label{t}
Let $N,p,q$  and $M$ be fixed. Then, there exist  $\widehat{S}_0=\widehat{S}_0(N,p,q,M)\in \er$ and $\varepsilon_0=\varepsilon_0(N,p,M)$, such that if $u $   is a solution of ({\ref{1}}) with blow-up graph $\Gamma:\{x\rightarrow T(x)\}$ and  $x_0$ is a non characteristic point, then

i)
 For all
 $s\ge \widehat{s}_0(x_0)=\max(\widehat{S}_0(N,p,q,M),-\log \frac{T(x_0)}4)$,
\begin{equation*}
0<\varepsilon_0\le \|w_{x_0,T(x_0)}(s)\|_{H^{1}(B)}+
\|\partial_s w_{x_0,T(x_0)}(s)\|_{L^{2}(B)}
\le K,
\end{equation*}
where $w_{x_0,T(x_0)}$ is defined in (\ref{B}) and $B$  is the unit ball of $\er^N$.\\
ii)  For all
  $t\in [t_0(x_0),T(x_0))$, where  $t_0(x_0)=\max (T(x_0)-e^{-\widehat{s}_0(x_0)},\frac{3T(x_0)}4)$, we have

\begin{eqnarray*}
&&0<\varepsilon_0\le (T(x_0)-t)^{\frac{2}{p-1}}\frac{\|u(t)\|_{L^2(B(x_0,{T(x_0)-t}))}}{ (T(x_0)-t)^{\frac{N}{2}}}\nonumber\\
&&+ (T(x_0)-t)^{\frac{2}{p-1}+1}\Big (\frac{\|u_t(t)\|_{L^2(B(x_0,{T(x_0)-t}))}}{ (T(x_0)-t)^{\frac{N}{2}}}+
 \frac{\|\grad u(t)\|_{L^2(B(x_0,{T(x_0)-t}))}}{ (T(x_0)-t)^{\frac{N}{2}}}\Big )\le K,
\end{eqnarray*}
where  $K=K(N,p,q,\widehat{s}_0(x_0),\|(u(t_0(x_0)),u_t(t_0(x_0)))\|_{
H^{1}\times L^{2}(B(x_0,\frac{e^{-\widehat{s}_0(x_0)}}{\delta_0(x_0)})
)})$ and\\
$\delta_0(x_0)\in (0,1)$ is defined in (\ref{cone1}).
\end{theo}
\no {\bf{Remark:}}

\no With this blow-up rate, one can ask whether the results proved by Merle and Zaag for the non perturbed problem in \cite{MZ4} \cite{MZ5} \cite{MZarxiv}, hold for equation (\ref{1}) (blow-up, profile, regularity of the blow-up graph, existence of characteristic points, etc...).  We believe that it is the case, however, the proof will be highly technical, with no interesting ideas (in particular, equation
(\ref{1}) is not conserved under the Lorentz transform, which is crucial in  \cite{MZ4} \cite{MZ5} \cite{MZarxiv}, and lots of minor term will
appear in the analysis). Once again, we believe that the key point in the analysis of blow-up for
equation (\ref{1}) is the derivation of a Lyapunov functional
in similarity variables, which is the object of our paper.

\medskip

\no
As in the particular case where $(f,g)\equiv (0,0)$, the proof of Theorem \ref{t} relies on four ideas (the existence of a Lyapunov functional,  interpolation in Sobolev spaces, some
Gagliardo-Nirenberg estimates and
a covering technique adapted to the geometric shape of the blow-up surface). It happens that adapting the proof of
\cite{MZ3} given in the non perturbed case (\ref{mu=0}) is straightforward, except for a key
argument, where we bound the $L^{p+1}$ space-time norm of $w$. Therefore, we only present that argument, and refer to
\cite{MZ1} and \cite{MZ3}  for the rest of the proof.

\medskip

\no
This paper is divided in two sections, each of them devoted to the proof of a Theorem.
\section{A Lyapunov functional for equation (\ref{B}) }
This section is divided in two parts:
  \begin{itemize}
\item  We first prove the existence of a Lyapunov functional for equation (\ref{B}).
 \item Then, we give a blow-up criterion for  equation (\ref{B}) based on the Lyapunov functional.
  \end{itemize}
Throughout this section, we consider $(w,w_s)\in
 {\cal
  C}([s_0,+\infty),{\cal H})$ where $w$ is a solution of  (\ref{B})  and   $s_0\in \er$. We aim at proving that the functional
  $H$ defined in (\ref{f1}) is a Lyapunov functional  for equation (\ref{B}), provided that $s\ge S_0$, for some $S_0=S_0(N,p,q,M)$.
We denote by $C$ a constant which depends
only on  $(p,q,N,M )$.
We denote the unit ball
of $\er^N$ by $B$.
\subsection{Existence of a Lyapunov functional }
\begin{lem}
\label{lyap0}
Let $N,p,q$  and $M$ be fixed. There exists  $S_1=S_1(N,p,q,M)\in \er$
such that, for all $ s_0\in \er$ and $w$ solution of
equation (\ref{B}) satisfying  $(w,w_s)\in
 {\cal
  C}([s_0,+\infty),{\cal H})$, we have
 the following inequality, for all $s\ge \max (s_0,S_1)$,
\begin{eqnarray}\label{lem21}
\frac{d}{ds}(E_0(w)+I (w))\!\!\!\!&\le&\!\!\!\!-\frac{3\alpha}{2}\int_{B}w_s^2\frac{\rho}{1-|y|^2}{\mathrm{d}}y
+\Sigma_ {0}(s),
\end{eqnarray}
where $\Sigma_0$ satisfies
\begin{eqnarray*}
%\label{sigm}
\Sigma_ {0}(s) &\le&\!\!\!\!C e^{-2s}\int_{B}|\grad w|^2(1-|y|^2)\rho{\mathrm{d}}y+C  e^{-2s}\int_{B}w^2\rho{\mathrm{d}}y\\
%\nonumber\\
&&+
Ce^{-\frac{2(p-q)s}{p-1}}\int_{B}\!\!|w|^{p+1}\rho {\mathrm{d}}y
 + C e^{-\frac{2(p-q)s}{p-1}}.\qquad \qquad
\end{eqnarray*}
\end{lem}
Proof:  Multipling $(\ref{B})$ by $w_s\rho$, and integrating over the ball $B$, we obtain, for all $s\ge s_0$,
(recall from \cite{AM} that in the case where, $(f,g)\equiv (0,0)$, we have\\ $\frac{d}{ds}E_0(w)=-2\alpha \int_{B}w_s^2\frac{\rho}{1-|y|^2}{\mathrm{d}}y$.)
\begin{eqnarray}\label{E0}
\frac{d}{ds}(E_0(w)+I (w))&=&-2\alpha \int_{B}w_s^2\frac{\rho}{1-|y|^2}{\mathrm{d}}y+
\underbrace{\frac{2(p+1)}{p-1}e^{-\frac{2(p+1)s}{p-1}}\int_{B}F\Big(e^{\frac{2s}{p-1}}w\Big)\rho {\mathrm{d}}y}_{I_1} \nonumber\\
&&+\underbrace{\frac{2}{p-1}e^{-\frac{2ps}{p-1}}\int_{B}f\Big(e^{\frac{2s}{p-1}}w\Big)w\rho {\mathrm{d}}y }_{I_2}\nonumber\\
&&+\underbrace{ e^{-\frac{2ps}{p-1}}\int_{B}g\Big(e^{\frac{(p+1)s}{p-1}}(w_s+y.\grad w+\frac{2}{p-1}w)\Big)w_s\rho {\mathrm{d}}y}_{I_3}.
\end{eqnarray}
By exploiting the fact that
$ |{F(x)}|+|x{f(x)}|\le C( 1+|x|^{q+1}),$ we obtain
\begin{eqnarray}\label{I10}
|I_1|+|I_2|&\le& Ce^{-\frac{2(p+1)s}{p-1}}\int_{B}(1+|e^{\frac{2s}{p-1}}w|^{q+1})\rho {\mathrm{d}}y\nonumber\\
&\le& Ce^{-\frac{2(p+1)s}{p-1}}+ Ce^{-\frac{2(p-q)s}{p-1}}\int_{B}|w|^{q+1}\rho {\mathrm{d}}y.
\end{eqnarray}
Noticing that $|x|^{q+1}\le C( 1+|x|^{p+1})$, we deduce from (\ref{I10}) that
for all $s \ge \max (s_0,0)$,
\begin{eqnarray}\label{I100}
|I_1|+|I_2|\le  Ce^{-\frac{2(p-q)s}{p-1}}+ Ce^{-\frac{2(p-q)s}{p-1}}\int_{B}|w|^{p+1}\rho {\mathrm{d}}y.
\end{eqnarray}
%on the one hand,
Since
$ |g(x)|\le M( 1+|x|)$, we write
\begin{eqnarray}\label{I122}
|I_3|&\le& Ce^{-s}
\int_{B}w_s^2\rho {\mathrm{d}}y +
Ce^{-s}\int_{B}|y.\grad w ||w_s|\rho {\mathrm{d}}y\nonumber\\
&&+ C e^{-s}\int_{B}|ww_s|\rho {\mathrm{d}}y+ C e^{-\frac{2ps}{p-1}}\int_{B}|w_s|\rho {\mathrm{d}}y.\qquad
\end{eqnarray}
By exploiting the inequality $ab\le \frac{a^2}{2}+\frac{b^2}{2},$ we obtain
 \begin{eqnarray}\label{I13}
C e^{-s}\int_{B}|y.\grad w| | w_s|\rho{\mathrm{d}}y\le
\frac{\alpha}{8}\int_{B}w_s^2\frac{\rho}{1-|y|^2}{\mathrm{d}}y
+C e^{-2s}\int_{B}|\grad w|^2(1-|y|^2)\rho{\mathrm{d}}y.
\end{eqnarray}
Similarly, we prove that
 \begin{equation}
 \label{I14}
 C e^{-s}\int_{B}|ww_s|\rho {\mathrm{d}}y\le \frac{\alpha}{8}\int_{B}w_s^2\frac{\rho}{1-|y|^2}{\mathrm{d}}y
+C  e^{-2s}\int_{B}w^2\rho{\mathrm{d}}y.
\end{equation}
We infer  from the inequality $ |a|\le 1+a^2$  that
\begin{equation}
 \label{I15}
 C e^{-\frac{2ps}{p-1}}\int_{B}|w_s|\rho {\mathrm{d}}y\le C  e^{-\frac{2ps}{p-1}}+C  e^{-\frac{2ps}{p-1}}\int_{B}w_s^2\frac{\rho}{1-y^2} {\mathrm{d}}y.
\end{equation}
Combining   (\ref{I122}), (\ref{I13}), (\ref{I14}) and (\ref{I15}), we conclude that,
for all $s \ge \max (s_0,0)$,
\begin{eqnarray}\label{I16}
|I_3|&\le& (Ce^{-s}+
\frac{\alpha}{4})\int_{B}w_s^2\frac{\rho}{1-|y|^2}{\mathrm{d}}y+C  e^{-2s}\int_{B}w^2\rho{\mathrm{d}}y\nonumber\\
&&+C e^{-2s}\int_{B}|\grad w|^2(1-|y|^2)\rho{\mathrm{d}}y+ C e^{-\frac{2ps}{p-1}}.
\end{eqnarray}
Then, by using  (\ref{E0}), (\ref{I100}) and (\ref{I16}), we deduce that,
for all $s \ge \max (s_0,0)$,
\begin{eqnarray}\label{E}
\frac{d}{ds}(E(w)+I (w))\!\!\!\!&\le&\!\!\!\!(-\frac{7\alpha}{4}+Ce^{-s}
)\int_{B}w_s^2\frac{\rho}{1-|y|^2}{\mathrm{d}}y
+C e^{-2s}\int_{B}|\grad w|^2(1-|y|^2)\rho{\mathrm{d}}y
\nonumber\\
&&+C  e^{-2s}\int_{B}w^2\rho{\mathrm{d}}y+   Ce^{-\frac{2(p-q)s}{p-1}}\int_{B}\!\!|w|^{p+1}\rho {\mathrm{d}}y+C e^{-\frac{2(p-q)s}{p-1}}.%\qquad\qquad
\end{eqnarray}
Taking $ S_1=S_1(N,p,q,M)$ large enough,  we have the estimate (\ref{lem21}). This concludes the proof of Lemma
 \ref{lyap0}.

\Box

\vspace{0.3cm}

We are now going to prove the following estimate for the functional $J$:
\begin{lem}
\label{lyap1}
Let $N,p,q$  and $M$ be fixed. There exists  $S_2=S_2(N,p,q,M)\in \er$
such that, for all $ s_0\in \er$ and $w$ solution of
equation (\ref{B}) satisfying  $(w,w_s)\in
 {\cal
  C}([s_0,+\infty),{\cal H})$,  $J$ satisfies
 the following inequality, for all $s\ge \max (s_0,S_2)$:
\begin{eqnarray}\label{theta1}
\frac{d}{ds}J(w)&\le &\frac{\alpha}{2}
\!\int_{B}\!w_s^2\frac{\rho}{1-|y|^2}{\mathrm{d}}y+\frac{p+3}2 e^{-\gamma s}E(s)\nonumber\\
&&-\frac{p-1}{4} e^{-\gamma s}\int_{B}(|\grad w|^2(1-|y|^2)\rho{\mathrm{d}}y\\
&&
- \frac{p+1}{2(p-1)}e^{-\gamma s}\int_{B}w^2\rho{\mathrm{d}}y
-\frac{p-1}{2(p+1)} e^{-\gamma s}\int_{B}|w|^{p+1}\rho{\mathrm{d}}y
+\Sigma_ {1}(s),\nonumber
\end{eqnarray}
where $\gamma =\min(\frac12,\frac{p-q}{p-1})>0$ and  $\Sigma_1(s)$ satisfies
\begin{eqnarray}
\label{sigm}
 \Sigma_ 1(s) &\le&
C e^{-2\gamma s}\int_{B}w^2\rho
{\mathrm{d}}y+
+C e^{-2\gamma s}\int_{B}\!|\grad w|^2\rho(1-|y|^2)
{\mathrm{d}}y
\nonumber\\
&&+C e^{-2\gamma s}\int_{B}| w|^{p+1}\rho{\mathrm{d}}y+ C e^{-2\gamma s}.
\end{eqnarray}
\end{lem}

Proof:
Note that  $J$ is a differentiable  function for  all
 $s\ge s_0$ and that
\begin{eqnarray*}
\frac{d}{ds}J(w)&=&\gamma  e^{-\gamma s}\int_{B}ww_{s}\rho{\mathrm{d}}y- e^{-\gamma s}\int_{B}w^2_{s}\rho{\mathrm{d}}y- e^{-\gamma s}\int_{B}
ww_{ss}\rho{\mathrm{d}}y.
\end{eqnarray*}
By using  equation (\ref{B}) and integrating by parts, we have
\begin{eqnarray}\label{t1}
\frac{d}{ds}J(w)&=&
- e^{-\gamma s}\int_{B}w^2_{s}\rho{\mathrm{d}}y
+ e^{-\gamma s}\int_{B}(|\grad w|^2-(y.\grad w)^2)\rho{\mathrm{d}}y\nonumber\\
&&+ \frac{2p+2}{(p-1)^2}e^{-\gamma s}\int_{B}w^2\rho{\mathrm{d}}y- e^{-\gamma s}\int_{B}|w|^{p+1}\rho{\mathrm{d}}y\nonumber\\
&&+
(\gamma +\frac{p+3}{p-1}-2N) e^{-\gamma  s}\int_{B}ww_s \rho{\mathrm{d}}y-2e^{-\gamma s}\int_{B}ww_s(y.\grad \rho){\mathrm{d}}y\nonumber\\
&&-2e^{-\gamma s}\int_{B}w_s(y.\grad w) \rho{\mathrm{d}}y-e^{-\frac{2ps}{p-1}-\gamma s}\int_{B}wf\Big(e^{\frac{2s}{p-1}}w\Big){\rho}{\mathrm{d}}y\\
&& -
e^{-\frac{2ps}{p-1}-\gamma s}\int_{B}wg\Big(e^{\frac{(p+1)s}{p-1}}(w_s+y.\grad w+\frac{2}{p-1}w)\Big){\rho}{\mathrm{d}}y.\nonumber
\end{eqnarray}
By combining (\ref{f}), (\ref{f2}) and (\ref{t1}),   we write
\begin{eqnarray}\label{theta}
\frac{d}{ds}J(w)&\le &\frac{p+3}2 e^{-\gamma s}E(s)
-\frac{p-1}{4} e^{-\gamma s}\int_{B}(|\grad w|^2-(y.\grad w)^2)\rho{\mathrm{d}}y\nonumber\\
&&- \frac{p+1}{2(p-1)}e^{-\gamma s}\int_{B}w^2\rho{\mathrm{d}}y
-\frac{p-1}{2(p+1)} e^{-\gamma s}\int_{B}|w|^{p+1}\rho{\mathrm{d}}y\nonumber\\
&&+
\underbrace{(\gamma +\frac{p+3}{p-1}-2N+\frac{p+3}2 e^{-\gamma s}) e^{-\gamma  s}\int_{B}ww_s \rho{\mathrm{d}}y}_{J_1}\nonumber\\
&&\underbrace{-2e^{-\gamma s}\int_{B}ww_s(y.\grad \rho){\mathrm{d}}y}_{J_2}
\underbrace{-2e^{-\gamma s}\int_{B}w_s(y.\grad w) \rho{\mathrm{d}}y}_{J_3}\nonumber\\
&&\underbrace{-e^{-\frac{2ps}{p-1}-\gamma s}\int_{B}wf\Big(e^{\frac{2s}{p-1}}w\Big){\rho}{\mathrm{d}}y}_{J_4}\nonumber\\
&&\underbrace{ -
e^{-\frac{2ps}{p-1}-\gamma s}\int_{B}wg\Big(e^{\frac{(p+1)s}{p-1}}(w_s+y.\grad w+\frac{2}{p-1}w)\Big){\rho}{\mathrm{d}}y}_{J_5}\nonumber\\
&&+
\underbrace{\frac{(p+3)}2  e^{-\frac{2(p+1)s}{p-1}-\gamma s}\displaystyle\int_{B}F(e^{\frac{2}{p-1}s}w)\rho {\mathrm{d}}y}_{J_6}.
\end{eqnarray}
We now study each of the last five terms. To estimate $J_1$, we use the fact that
 for  all
 $s\ge \max(s_0,0)$,
 \begin{equation}\label{es}
\big |\gamma +\frac{p+3}{p-1}-2N+\frac{p+3}2 e^{-\gamma s}\big |\le C.
\end{equation}
By using (\ref{es}) and the Cauchy-Schwartz inequality we obtain
\begin{eqnarray}\label{J1}
|J_1|& \le& C  e^{-\gamma s}
\int_{B}|ww_s| \rho{\mathrm{d}}y
\le
 \frac{\alpha}{8}\int_{B}w_s^2\frac{\rho}{1-|y|^2}{\mathrm{d}}y+
C e^{-2\gamma s}\int_{B}w^2\rho
{\mathrm{d}}y.\qquad\qquad
\end{eqnarray}
Now we estimate the expression $J_2$.
Since we have $\ds{y.\grad \rho=-2\alpha \frac{|y|^2}{(1-|y|^2)}\rho,}$
we can use the Cauchy-Schwartz inequality to write
\begin{eqnarray}\label{J2}
|J_2|& \le& C e^{-\gamma s}
\int_{B}|w_s|(1-|y|^2)^{\frac{\alpha-1}{2}}|w||y|(1-|y|^2)^{\frac
{\alpha-1}{2}}{\mathrm{d}}y\nonumber\\
&\le&\frac{\alpha}{8}\int_{B}w_s^2\frac{\rho}{1-|y|^2}{\mathrm{d}}y+C e^{-2\gamma s}\int_{B}w^2\frac{|y|^2\rho}{1-|y|^2}
{\mathrm{d}}y.
\end{eqnarray}
Since we have the following Hardy type inequality for any $w\in
H^{1}_{loc,u}(\er^N)$ (see appendix $B$  in \cite{MZ1} for details):
\begin{eqnarray}\label{hardyJJ}
\int_{B}w^2\frac{|y|^2\rho}{1-|y|^2}{\mathrm{d}}y
&\le&C\int_{B}|\grad w|^2\rho(1-|y|^2){\mathrm{d}}y+C \int_{B}w^2\rho
{\mathrm{d}}y,
\end{eqnarray}
we use  (\ref{J2}) and (\ref{hardyJJ}) to conclude that
\begin{equation}\label{J22}
|J_2|\le\frac{\alpha}{8}\!\int_{B}\!w_s^2\frac{\rho}{1-|y|^2}{\mathrm{d}}y+C
e^{-2\gamma s}\int_{B}\!w^2\rho
{\mathrm{d}}y
+C e^{-2\gamma s}\int_{B}\!|\grad w|^2\rho(1-|y|^2)
{\mathrm{d}}y.\quad
\end{equation}
By using the Cauchy-Schwartz inequality, we have
\begin{eqnarray}\label{J3}
|J_3|&\le&\frac{\alpha}{8}\int_
{B}w_s^2\frac{\rho}{1-|y|^2}{\mathrm{d}}y
+C e^{-2\gamma s}\int_{B}|\grad w|^2\rho(1-|y|^2){\mathrm{d}}y.
\end{eqnarray}
By exploiting the fact that
$ |{F(x)}|\le CM( 1+|x|^{q+1})$  and
$ |{f(x)}|\le M(1+|x|^q),$
we write
\begin{eqnarray}\label{J4}
|J_4|+|J_6|&\le&
C e^{-\frac{2(p-q)s}{p-1}-\gamma s}\int_{B}(1+| w|^{q+1})\rho{\mathrm{d}}y
\le C e^{-2\gamma s}\int_{B}(1+|w|^{p+1})\rho{\mathrm{d}}y\nonumber\\
&\le&C e^{-2\gamma s}+
 C e^{-2\gamma s}\int_{B}|w|^{p+1}\rho{\mathrm{d}}y.
\end{eqnarray}
In a similar way,  by  using the fact that
$ |g(x)|\le M(1+|x|)$, we write
\begin{eqnarray*}\label{J5}
|J_5|&\le& Ce^{-2\gamma s}
\int_{B}w_s^2\rho {\mathrm{d}}y +
Ce^{-2\gamma s}\int_{B}|y.\grad w ||w|\rho {\mathrm{d}}y\nonumber\\
&&+ C e^{-2\gamma s}\int_{B}w^2\rho {\mathrm{d}}y+ C e^{-2\gamma s}.\qquad
\end{eqnarray*}
Then, by (\ref{hardyJJ}), we  have
\begin{eqnarray}\label{J55}
|J_5|&\le& Ce^{-2\gamma s}
\int_{B}w_s^2\frac{\rho}{1-|y|^2} {\mathrm{d}}y +C e^{-2\gamma s}\int_{B}|\grad w|^2\rho(1-|y|^2)
{\mathrm{d}}y\nonumber\\
&& +Ce^{-2\gamma s}\int_{B}w^2\rho {\mathrm{d}}y+ C e^{-2\gamma s}.
\end{eqnarray}
Finally, by using (\ref{theta}), (\ref{J1}), (\ref{J22}), (\ref{J3}), (\ref{J4}) and (\ref{J5}) we deduce that
\begin{eqnarray*}%\label{theta11}
\frac{d}{ds}J(w)&\le &\frac{p+3}2 e^{-\gamma s}E(s)
-\frac{p-1}{4} e^{-\gamma s}\int_{B}(|\grad w|^2-(y.\grad w)^2)\rho{\mathrm{d}}y\nonumber\\
&&- \frac{p+1}{2(p-1)}e^{-\gamma s}\int_{B}w^2\rho{\mathrm{d}}y
-\frac{p-1}{2(p+1)} e^{-\gamma s}\int_{B}|w|^{p+1}\rho{\mathrm{d}}y\nonumber\\
&&+
C e^{-2\gamma s}\int_{B}w^2\rho
{\mathrm{d}}y+
+C e^{-2\gamma s}\int_{B}\!|\grad w|^2\rho(1-|y|^2)
{\mathrm{d}}y\nonumber\\
&&+(
\frac{3\alpha}{8}+Ce^{-2\gamma s})\!\int_{B}\!w_s^2\frac{\rho}{1-|y|^2}{\mathrm{d}}y+
 C e^{-2\gamma s}+C e^{-2\gamma s}\int_{B}| w|^{p+1}\rho{\mathrm{d}}y.
\end{eqnarray*}
Since $|y.\grad w|\le |y||\grad w|$, it follows that
\begin{eqnarray}\label{last}
 \int_{B}|\grad w|^2\rho(1-|y|^2){\mathrm{d}}y&\le &\int_{B}(|\grad w|^2-(y.\grad w)^2) \rho{\mathrm{d}}y.
\end{eqnarray}
Taking $S_2=S_2(N,p,q,M)$ large enough,  we have easily  the estimate (\ref{theta1}) and (\ref{sigm}).
This concludes the proof of  Lemma \ref{lyap1}.\Box

\bigskip

With Lemmas \ref{lyap0} and \ref{lyap1}, we are in a position to prove Theorem \ref{lyap}.

\bigskip

\no {\bf{Proof of Theorem \ref{lyap}}}\\
From Lemmas \ref{lyap0} and \ref{lyap1}, we obtain for all $s\ge \max(s_0, S_1, S_2)$,
%Take $S_0\ge \max(S_1,S_2)$, where $S_1$ and $S_2$ are defined in
%Lemma \ref{lyap} and  Lemma \ref{lyap1}, we obtain
\begin{eqnarray}\label{E11}
\frac{d}{ds}E(w)&\le &C e^{-2\gamma s}+\frac{p+3}2 e^{-\gamma s}E(w)
-\alpha\!\int_{B}\!w_s^2\frac{\rho}{1-|y|^2}{\mathrm{d}}y\nonumber\\
&&+(C e^{-\gamma s}-\frac{p-1}{4}) e^{-\gamma s}\int_{B}(|\grad w|^2(1-|y|^2)\rho{\mathrm{d}}y\nonumber\\
&&+(C e^{-\gamma s}
- \frac{p+1}{2(p-1)})e^{-\gamma s}\int_{B}w^2\rho{\mathrm{d}}y\nonumber\\
&&+(C e^{-\gamma s}-\frac{p-1}{2(p+1)}) e^{-\gamma s}\int_{B}|w|^{p+1}\rho{\mathrm{d}}y.
%,\quad \forall s\ge s_0.
\end{eqnarray}
We now choose $S_0\ge \max(S_1,S_2)$, large enough, so that for all $s\ge S_0$, we have
\begin{eqnarray*}
 \frac{p-1}{4}-C e^{-\gamma s} \ge 0, \qquad
 \frac{p+1}{2(p-1)}
-C e^{-\gamma s} \ge 0,\quad \frac{p-1}{2(p+1)}-C e^{-\gamma s}\ge 0.
\end{eqnarray*}
Then, we deduce that,
 for all $s\ge \max(S_0,s_0)$,
 we have
\begin{eqnarray}\label{E111}
\frac{d}{ds}E(w)&\le & C e^{-2\gamma s}+\frac{p+3}2 e^{-\gamma s}E(w)-\alpha \int_{B}w_s^2\frac{\rho}{1-|y|^2}{\mathrm{d}}y.
\end{eqnarray}
Finally, we prove easily that the
function $H$ satisfies, for all $s\ge \max(S_0,s_0)$,
\begin{eqnarray}\label{H01}
\frac{d}{ds}H(w)&\le & (Ce^{\frac{ p+3}{2\gamma} e^{-\gamma s}}-2\theta \gamma ) e^{-2\gamma s}-\alpha e^{\frac{ p+3}{2\gamma} e^{-\gamma s}}\int_{B}w_s^2\frac{\rho}{1-|y|^2}{\mathrm{d}}y\nonumber\\
&\le & (C-2\theta \gamma ) e^{-2\gamma s}-\alpha \int_{B}w_s^2\frac{\rho}{1-|y|^2}{\mathrm{d}}y.
\end{eqnarray}
We now choose $\theta$ large enough,  so we have $C-2\theta \gamma \le 0$ and then
\begin{eqnarray}\label{H02}
\frac{d}{ds}H(w)
&\le & -\alpha \int_{B}w_s^2\frac{\rho}{1-|y|^2}{\mathrm{d}}y.
\end{eqnarray}
Now (\ref{L1}) is a direct consequence of inequality
(\ref{H02}).\\
This concludes the proof of Theorem \ref{lyap}.

\Box
\subsection{ A blow-up criterion  in the $w(y,s)$ variable}

\no We now claim the following proposition:
\begin{lem}\label{L02}
Let $N,p,q,M$  be fixed. There exists $S_3=S_3(N,p,q,M)\ge S_0$ such that,
for all $ s_0\in \er$ and $w$ solution of
equation (\ref{B}) defined to the left of $s_0$, such that  $\|w(s)\|_{L^{p+1}(B)}$ is locally bounded, if
  $H(w(s_3))<0$ for some $s_3\ge \max(S_3,s_0)$, then $w$
blows up in some finite time
$S>s_3$.
\end{lem}

\no {\bf{Remark}}\\
If $w=w_{x_0,T_0}$ defined from a solution of (\ref{1}) by (\ref{scaling})
and $x_0$ is non characteristic point, then
$\|w(s)\|_{H^1(B)}$ is locally bounded and  so  is $\|w(s)\|_{L^{p+1}(B)}$ by
Sobolev's embedding.

\no Proof:
The argument is the same as in the corresponding part
in \cite{AM}. We write the proof for completeness. Arguing by contradiction,
we assume that there exists a
solution $w$ on $B$, defined for all time $s\in [s_3,+\infty[$, where  $H(w(s_3))<0$.
Since the energy $H$ decreases in time, we have  $H(w(1+s_3))<0$.

\no
Consider now for $\delta>0$ the function $\widetilde{w}^{\delta}(y,s)$
for $(y,s)\in  B\times [1+s_3,+\infty[$ defined by
\begin{equation*}
\forall s\ge 1+s_3, \forall y\in B,
\widetilde{w}^{\delta}(y,s)=\frac{1}{(1+\delta e^s)^{\frac2{p-1}}} w(
\frac{y}{1+\delta e^s},-\log(\delta+e^{-s})).
\end{equation*}
\begin{itemize}
\item  (A) Note that $\widetilde{w}^{\delta}$ is defined in  $ B\times [1+s_3,+\infty[$,
whenever $\delta>0$ is small enough such that  $-\log(\delta+e^{-1-s_3})\ge s_3.$
\item (B) From  its construction, $\widetilde{w}^{\delta}$
is also a solution of (\ref{B}) (Indeed, let $u$ be such that $w=w_{0,0}$ in defintion  (\ref{scaling}). Then $u$ is a solution  of (\ref{1}) and $\widetilde{w}^{\delta}=w_{0,-\delta}$ is defined as in (\ref{scaling}); so
$\widetilde{w}^{\delta}$
is  a solution of (\ref{B})).
\item (C) For $\delta$ small enough, we have
$H(\widetilde{w}^{\delta}(1+s_3))<0$ by continuity of the function
$\delta \mapsto H(\widetilde{w}^{\delta}(1+s_3))$.
Then, we write that
$H(\widetilde{w}^{\delta}(1+s_3))<0$.
\end{itemize}
Now, we
fix $\delta=\delta_0>0$ such that (A), (B) and (C) hold.
Let us note that we have
\begin{eqnarray}\label{c1}
- e^{-\gamma s}\int_{B}w^{\delta_0}w^{\delta_0}_s\rho {\mathrm{d}}y&\ge&-\frac{1}{4}
\int_{B}(w^{\delta_0}_s)^2\rho{\mathrm{d}}y- e^{-2\gamma s}
\int_{B}(w^{\delta_0})^2\rho{\mathrm{d}}y
\end{eqnarray}
and from (\ref{I100})
\begin{eqnarray}\label{c2}
- e^{-\frac{2(p+1)s}{p-1}}\displaystyle\int_{B}F(e^{\frac{2}{p-1}s}w^{\delta_0})\rho {\mathrm{d}}y&\ge&
-C e^{-2\gamma s}-C e^{-2\gamma s}\int_{B}|w^{\delta_0}|^{p+1}\rho{\mathrm{d}}y.
\end{eqnarray}
By (\ref{f}), (\ref{f2}), (\ref{c1}) and (\ref{c2}) we deduce
\begin{eqnarray}
E(w^{\delta_0}(s))&\ge& \frac{1}{4}
\int_{B} (w_s^{\delta_0})^2\rho {\mathrm{d}}y
+(\frac{p+1}{(p-1)^2}- e^{-2\gamma s})\int_{B}(w^{\delta_0})^2\rho {\mathrm{d}}y\nonumber\\
&&-(\frac{1}{p+1}+C e^{-2\gamma s})\int_{B}|w^{\delta_0}|^{p+1} \rho {\mathrm{d}}y
-C e^{-2\gamma s}.
\end{eqnarray}
We now choose $s_4\ge s_3$ large enough, so that we have
$\frac{p+1}{(p-1)^2}- e^{-2\gamma s_4} \ge 0.$
Then, we deduce that we have, for all $s\ge s_4$,
\begin{equation*}
E(w^{\delta_0}(s))\ge
-(\frac{1}{p+1}+C e^{-2\gamma s})\int_{B}|w^{\delta_0}|^{p+1} \rho {\mathrm{d}}y
-C e^{-2\gamma s}.
\end{equation*}
Since $\rho\le 1$, after a change of variables, we find that
\[
E(w^{\delta_0}(s))
\ge-\frac{(\frac{1}{p+1}+C e^{-2\gamma s})}
{(1+\delta_0e^s)^{\frac{4}{p-1}+2-N}}\displaystyle\int_{B}|w(z,-\log (\delta_0+e^{-s}))|^{p+1} {\mathrm{d}}z-C e^{-2\gamma s}.
\]
Since we have $-\log (\delta_0+e^{-s})\rightarrow -\log (\delta_0)$ as $s\rightarrow +\infty$ and since $\|w(s)\|_{L^{p+1}(B)}$ is locally bounded by hypothesis, by a continuity argument,  it follows that the former integral remains bounded and
\begin{eqnarray*}%\label{}
E(w^{\delta_0}(s))&\ge&
-\frac{C}{(1+\delta_0e^s)^{\frac{4}{p-1}+2-N}}-C e^{-2\gamma s}\rightarrow 0,
\end{eqnarray*}
as $s\rightarrow +\infty$ (use the fact that $\frac{4}{p-1}+2-N>0$ which follows from the fact that $p<p_c$).
So, from (\ref{f1}), it follows that
\begin{eqnarray}\label{d1}
\liminf_{s\rightarrow +\infty}H(w^{\delta_0}(s))\ge 0 .
\end{eqnarray}
\no The  inequality (\ref{d1}) contradicts the inequality  $H(w^{\delta_0}(s_3+1))<0$ and the fact that the
 energy $H$ decreases in time for $s\ge s_3$. This concludes the proof of Lemma \ref{L02}.

\Box

\section{Boundedness of the solution in similarity variables}
\no We prove Theorem 2.2 here. Note that the lower bound follows from the finite speed of propagation and wellposedness  in $H^1\times L^2$.
For a detailed argument in the similar case of equation (\ref{mu=0}), see Lemma 3.1 (page 1136) in \cite{MZ3}.\\
\no We consider $u$  a solution of (\ref{1}) which is defined under the graph of $x\mapsto T(x)$, and $x_0$ a non characteristic point.
Given some $T_0\in (0,T(x_0)]$, we introduce $w_{x_0,T_0}$ defined in (\ref{scaling}), and write $w$ for simplicity, when there is no ambiguity.
We aim   at bounding $\|(w,\partial_s w)(s)\|_{H^1\times L^2(B)}$ for $s$ large.

%\no {\bf{Remarque:}} The function $w$ is a global solution of equation (\ref{B}). Moreover, there exist $M_0>0$, such %that we have
%$$\|(u(s),u_{\tau}(s))\|_{H^{1}_{loc,u}(\er^N)\times L^{2}_{loc,u}(\er^N)}\le M_0,\qquad \forall s\in [-\log T,s_0] %$$

\no As in \cite{MZ1}, by combining Theorem \ref{lyap} and  Lemma \ref{L02} (use in particular the remark after that Lemma) we get the following bounds:
\begin{cor} (Bounds on $E$) For all $s\ge \widehat{s}_3=\widehat{s}_3(T_0)=\max(S_3,-\log T_0)$, $s_2\ge s_1\ge \widehat{s}_3$, it holds that
\label{cor1}
\begin{eqnarray}\label{cor01}
-C\le E(w(s))\le  M_0\nonumber\\
\int_{s_1}^{s_2}\int_{B}w_s^2(y,s)\frac{\rho}{1-|y|^2}dy ds\le M_0,
\end{eqnarray}
where  $M_0=M_0(N,p,q,M,\widehat{s}_3(T_0),\|(u(t_3),u_t(t_3))\|_{
H^{1}\times L^{2}(B(x_0,\frac{e^{-\widehat{s}_3(T_0)}}{\delta_0(x_0)})
)})$,\\
 $t_3=t_3(T_0)=T_0-e^{-\widehat{s}_3(T_0)}$,
$C=C(N,p,q,M)$    and   $\delta_0(x_0)\in (0,1)$ is defined in (\ref{cone1}).
\end{cor}
\no Starting from these bounds, the proof of Theorem 1.2 is similar to the proof in \cite{MZ3}
except for
the treatment
of the perturbation terms. In our opinion, handling these terms is straightforward in
all the steps of the proof, except for the first step, where we bound the time averages of the $L^{p+1}_{\rho}(B)$ norm of $w$. For that reason,
we only give that step and refer to  \cite{MZ3} for the remaining steps in the proof of  Theorem 1.2. This is the step we prove here
(In the following
$K_1$ denotes  a constant that depends only on $p$, $q$, $N$, $M$,  $C$, $M_0$, and $\varepsilon$ is an arbitrary positive number in $]0,1[$).
\begin{pro}\label{pro}(Control of the space-time $L^{p+1}$ norm of $w$)\\
 For all  $s\ge 1+\widehat{s}_3$,
\begin{equation}\label{pro1}
\int_{s}^{s+1}\int_{B}\!\!|w|^{p+1}{\rho}{\mathrm{d}}y{\mathrm{d}}s\le K_1(M_0,C,N,p,q,M ).
\end{equation}
\end{pro}
Proof:
For $s\ge 1+\widehat{s}_3$, let us work with time integrals betwen $s_1$ et $s_2$ where $s_1\in [s-1,s]$
and $s_2\in [s+1,s+2]$.
 By integrating the expression (\ref{f}) of $E$ in time between $s_1$ and $s_2$, where $s_2>s_1>\widehat{s}_3$, we obtain:
\begin{eqnarray}\label{et}
\int_{s_1}^{s_2}\!\!E(s)  ds&=&\displaystyle\int_{s_1}^{s_2}\int_{B}\!\!\Big ( \frac{1}{2}w_s^2
+\frac{p+1}{(p-1)^2}w^2
-\frac{1}{p+1}|w|^{p+1}\Big )\rho {\mathrm{d}}y{\mathrm{d}}s\\
&&+\frac{1}{2}\displaystyle\int_{s_1}^{s_2}\!\!\int_{B}\!\!\Big (|\grad w|^2-(y.\grad w)^2\Big )\rho {\mathrm{d}}y{\mathrm{d}}s-\int_{s_1}^{s_2}\!\! e^{-\gamma s}\!\!\displaystyle\int_{B}\!\!ww_s\rho {\mathrm{d}}y{\mathrm{d}}s\nonumber\\
&&-\int_{s_1}^{s_2}\!\! e^{-\frac{2(p+1)s}{p-1}}\displaystyle\int_{B}F(e^{\frac{2}{p-1}s}w)\rho {\mathrm{d}}y\mathrm{d}s.\nonumber
\end{eqnarray}
By multiplying the equation (\ref{B}) by $w\rho$ and integrating both in time and in space over $B\times [s_1,s_2]$, we obtain the following identity, after some
integration by parts :
\begin{eqnarray}\label{et1}
&&\Big [\int_{B}\!\!\Big (ww_{s}+(\frac{p+3}{2(p-1)}-N)w^2\Big ) \rho{\mathrm{d}}y\Big ]_{s_1}^{s_2}=
\int_{s_1}^{s_2}\!\!\int_{B}\!\!w^2_{s}\rho{\mathrm{d}}y{\mathrm{d}}s\nonumber\\
&&-\int_{s_1}^{s_2}\!\!\int_{B}\!\!(|\grad w|^2-(y.\grad w)^2)\rho{\mathrm{d}}y{\mathrm{d}}s
-\frac{2p+2}{(p-1)^2}\int_{s_1}^{s_2}\!\!\int_{B}\!\!w^2\rho{\mathrm{d}}y{\mathrm{d}}s\nonumber\\
&&+
\int_{s_1}^{s_2}\!\!\int_{B}\!\!|w|^{p+1}\rho{\mathrm{d}}y{\mathrm{d}}s+\!2\!\int_{s_1}^{s_2}\!\!\int_{B}\!\!ww_s(y.\grad \rho){\mathrm{d}}y{\mathrm{d}}s
+2\!\!\int_{s_1}^{s_2}\!\!\!\int_{B}\!\!\!w_s(y.\grad w) \rho{\mathrm{d}}y{\mathrm{d}}s\nonumber\\
&&+\!\int_{s_1}^{s_2}\!\!\int_{B}\!\!e^{-\frac{2ps}{p-1}}f\Big(e^{\frac{2s}{p-1}}w\Big)w\rho{\mathrm{d}}y{\mathrm{d}}s\nonumber\\
&& +\!\int_{s_1}^{s_2}\!\!\int_{B}\!\!
e^{-\frac{2ps}{p-1}}g\Big(e^{\frac{(p+1)s}{p-1}}(w_s+y.\grad w+\frac{2}{p-1}w)\Big)w\rho{\mathrm{d}}y{\mathrm{d}}s.\qquad
\end{eqnarray}
By combining the identities (\ref{et}) and (\ref{et1}), we obtain
\begin{eqnarray}\label{controlp}
&&\frac{(p-1)}{2(p+1)}\int_{s_1}^{s_2}\!\!\int_{B}\!\!|w|^{p+1}\rho{\mathrm{d}}y{\mathrm{d}}s\\
&=&\frac{1}{2}\Big [\int_{B}\!\!\Big (ww_{s}+(\frac{p+3}{2(p-1)}-N)w^2\Big ) \rho{\mathrm{d}}y\Big ]_{s_1}^{s_2}-
\int_{s_1}^{s_2}\!\!\int_{B}\!\!w^2_{s}\rho{\mathrm{d}}y{\mathrm{d}}s\nonumber\\
&&+\int_{s_1}^{s_2}\!\!E(s) ds
-\int_{s_1}^{s_2}\!\!\int_{B}\!\!ww_s(y.\grad \rho){\mathrm{d}}y{\mathrm{d}}s
-\int_{s_1}^{s_2}\!\!\int_{B}\!\!w_s(y.\grad w) \rho{\mathrm{d}}y{\mathrm{d}}s\nonumber\\
&& -\underbrace{\frac12\!\int_{s_1}^{s_2}\!\!\int_{B}\!\!
e^{-\frac{2ps}{p-1}}g\Big(e^{\frac{(p+1)s}{p-1}}(w_s+y.\grad w+\frac{2}{p-1}w)\Big)w\rho{\mathrm{d}}y{\mathrm{d}}s}_{A_1}\nonumber\\
&&-\underbrace{\int_{s_1}^{s_2}\!\! e^{-\gamma s}\!\!\displaystyle\int_{B}\!\!ww_s\rho {\mathrm{d}}y{\mathrm{d}}s}_{A2}
-\underbrace{\frac12\!\int_{s_1}^{s_2}\!\!\int_{B}\!\!e^{-\frac{2ps}{p-1}}f\Big(e^{\frac{2s}{p-1}}w\Big)w\rho{\mathrm{d}}y{\mathrm{d}}s}_{A_3}
\nonumber\\
&&+\underbrace{\int_{s_1}^{s_2}\!\! e^{-\frac{2(p+1)s}{p-1}}\displaystyle\int_{B}F(e^{\frac{2}{p-1}s}w)\rho {\mathrm{d}}y\mathrm{d}s}_{A4}.
\end{eqnarray}
We claim that Proposition \ref{pro} follows from the following Lemma where we    control all the terms on the right-hand side of the relation
 (\ref{controlp}) in terms of
the space-time $L^{p+1}$ norm of $w$:

\begin{lem}\label{g}
 For all  $s\ge 1+\widehat{s}_4$,
for some $\widehat{s}_4\ge \widehat{s}_3$, for all $\varepsilon>0$,
\begin{equation}\label{control}
\int_{s_1}^{s_2}\int_{B}\!\!|\grad w|^2(1-|y|^2)\rho{\mathrm{d}}y{\mathrm{d}}s\le K_1 +C
\int_{s_1}^{s_2}\!\!\int_{B}\!\!| w|^{p+1}\rho{\mathrm{d}}y{\mathrm{d}}s,
\end{equation}
\begin{equation}\label{control1}
\sup_{s\in [s_1,s_2]}\int_{B} \!\!w^2(y,s)\rho{\mathrm{d}}y\le \frac{K_1}{\varepsilon} +K_1\varepsilon
\int_{s_1}^{s_2}\!\!\int_{B}\!\!|w|^{p+1}\rho{\mathrm{d}}y{\mathrm{d}}s.
\end{equation}
\begin{equation}\label{control3}
\int_{s_1}^{s_2}\!\!\!\int_{B}\!|w_s y.\grad w| \rho{\mathrm{d}}y{\mathrm{d}}s
\le \frac{K_1}{\varepsilon} +K_1\varepsilon
\int_{s_1}^{s_2}\!\!\!\int_{B}\!| w|^{p+1}\rho{\mathrm{d}}y{\mathrm{d}}s,\qquad
\end{equation}
\begin{equation}\label{control30}
\int_{s_1}^{s_2}\!\!\!\int_{B}\!|w_sw y.\grad \rho|{\mathrm{d}}y{\mathrm{d}}s
\le \frac{K_1}{\varepsilon} +K_1\varepsilon
\int_{s_1}^{s_2}\!\!\!\int_{B}\!| w|^{p+1}\rho{\mathrm{d}}y{\mathrm{d}}s,\qquad
\end{equation}
\begin{equation}\label{control4}
\int_{B}|ww_s|\rho{\mathrm{d}}y\le
\int_{B} w_s^2\rho{\mathrm{d}}y+\frac{K_1}{\varepsilon} +K_1\varepsilon
\int_{s_1}^{s_2}\!\!\!\int_{B}\!| w|^{p+1}\rho{\mathrm{d}}y{\mathrm{d}}s,
\end{equation}
\begin{equation}\label{control5}
\int_{B}( w_s^2(y,s_1)+w_s^2(y,s_2))\rho{\mathrm{d}}y
\le K_1,
\end{equation}
\begin{eqnarray}\label{A10}
|A_1|&\le& \frac{K_1}{\varepsilon} +(K_1\varepsilon+
C  e^{-s_1})
\int_{s_1}^{s_2}\!\!\!\int_{B}\!| w|^{p+1}\rho{\mathrm{d}}y{\mathrm{d}}s,
\end{eqnarray}
\begin{eqnarray}\label{A20}
|A_2|
&\le &
\frac{K_1}{\varepsilon} +K_1\varepsilon
\int_{s_1}^{s_2}\!\!\!\int_{B}\!| w|^{p+1}\rho{\mathrm{d}}y{\mathrm{d}}s,\qquad
\end{eqnarray}
\begin{eqnarray}\label{A30}
|A_3|+|A_4|
&\le&  C+ Ce^{-\gamma s_1}\int_{s_1}^{s_2} \!\!\int_{B}|w|^{p+1}\rho {\mathrm{d}}y{\mathrm{d}}s.
\end{eqnarray}
\end{lem}
Indeed, from (\ref{controlp}) and this Lemma,
we deduce that
\begin{eqnarray*}%\label{}
\int_{s_1}^{s_2}\!\!\!\int_{B}\!| w|^{p+1}\rho{\mathrm{d}}y{\mathrm{d}}s&\le& \frac{K_1}{\varepsilon}+
(K_1\varepsilon+
C  e^{-s_1}+C  e^{-\gamma s_1})
\int_{s_1}^{s_2}\!\!\!\int_{B}\!| w|^{p+1}\rho{\mathrm{d}}y{\mathrm{d}}s.\qquad
\end{eqnarray*}
Taking $\widehat{s}_5$ large enough and $\varepsilon$ small enough so that
% for all $s_1\ge \widehat{s}_5$,
$C e^{-\widehat{s}_5}+C  e^{-\gamma \widehat{s}_5}\le \frac14$ and $K_1 \varepsilon \le \frac14$, we obtain (\ref{pro1}).

\bigskip

\no
It remains to prove Lemma 3.3.
%This concludes the proof of Proposition \ref{pro}.\\

\bigskip

\no
 Proof of Lemma 3.3:
 For the estimates (\ref{control}), (\ref{control1}), (\ref{control3}), (\ref{control30}), (\ref{control4}) and (\ref{control5}),  we can adapt with no difficulty
 the proof given in the case of  the wave equation treated in \cite{MZ1}.

\no
Now, we control the terms $A_1$, $A_2$, $A_3$ and $A_4$. Since
$ |g(x)|\le M(1+|x|)$, we write
\begin{eqnarray}
\label{I12}
|A_1|&\le& C\int_{s_1}^{s_2}\!\!e^{-s}
\int_{B}w_s^2\rho {\mathrm{d}}y\mathrm{d}s +
C\int_{s_1}^{s_2}\!\!e^{-s}\int_{B}|y.\grad w ||w|\rho {\mathrm{d}}y\mathrm{d}s\nonumber\\
&&+ C \int_{s_1}^{s_2}\!\!e^{-s}\int_{B}w^2\rho {\mathrm{d}}y\mathrm{d}s+ C \int_{s_1}^{s_2}\!\!e^{-\frac{2ps}{p-1}}\int_{B}|w|\rho {\mathrm{d}}y\mathrm{d}s.\qquad
\end{eqnarray}
By using (\ref{cor01}), we write
\begin{eqnarray}\label{CC0}
 C\int_{s_1}^{s_2}\!\!e^{-s}
\int_{B}w_s^2\rho {\mathrm{d}}y\mathrm{d}s \le K_1
\end{eqnarray}
Using the fact that $ e^{-s_1}\le1$ and the inequality (\ref{control1}), we obtain,
\begin{eqnarray}
\label{CC1}
C\!\!\int_{s_1}^{s_2}\!\!e^{-s}\!\!\int_{B}\!\!w^2{\rho}{\mathrm{d}}y{\mathrm{d}}s \le
C\sup_{s\in [s_1,s_2]}\int_{B}w^2{\rho}{\mathrm{d}}y
\le \frac{K_1}{\varepsilon} +K_1\varepsilon
\int_{s_1}^{s_2}\!\!\!\int_{B}\!| w|^{p+1}\rho{\mathrm{d}}y{\mathrm{d}}s.
%\qquad\qquad
\end{eqnarray}
We infer  from (\ref{CC1}) and the inequality $ |a|\le 1+a^2$  that
\begin{eqnarray}
 \label{CC2}
C  \int_{s_1}^{s_2}e^{-\frac{2ps}{p-1}}\int_{B}|w|\rho {\mathrm{d}}y\mathrm{d}s&\le& C  +C\int_{s_1}^{s_2}  e^{-s}\int_{B}w^2\rho {\mathrm{d}}y\mathrm{d}s\nonumber\\
&\le&\frac{K_1}{\varepsilon} +K_1\varepsilon
\int_{s_1}^{s_2}\!\!\int_{B}\!\!|w|^{p+1}\rho{\mathrm{d}}y{\mathrm{d}}s.
\end{eqnarray}
%On the other hand, by
Using the Cauchy-Schwarz inequality, we write
\begin{eqnarray}\label{control6}
C\!\! \int_{s_1}^{s_2}\!\!e^{-s}\!\!\int_{B}\!\!|w||y.\grad w|{\rho}{\mathrm{d}}y{\mathrm{d}}s\le
C\!\! \int_{s_1}^{s_2}\!\! e^{-s}\int_{B}\!\!|w||y||\grad w|{\rho}{\mathrm{d}}y{\mathrm{d}}s&&\nonumber\\
\le C \!\!\int_{s_1}^{s_2}\!\! e^{-s}\int_{B}\!\!w^2\frac{|y|^2}{1-|y|^2}{\rho}{\mathrm{d}}y{\mathrm{d}}s
+C e^{-s_1}\!\!\int_{s_1}^{s_2}\!\!\int_{B}\!\!|\grad w|^2(1-|y|^2){\rho}{\mathrm{d}}y{\mathrm{d}}s.&&
\end{eqnarray}
%Let's recall the Hardy  inequality:
%\begin{eqnarray}\label{control7}
%\int_{B}w^2\frac{|y|^2\rho}{1-|y|^2}{\mathrm{d}}y
%&\le&C\int_{B}|\grad w|^2\rho(1-|y|^2){\mathrm{d}}y+C \int_{B}w^2\rho
%{\mathrm{d}}y.
%\end{eqnarray}
By combining (\ref{CC1}),
(\ref{control6}),
%(\ref{control7})
(\ref{hardyJJ}),
(\ref{control}) and (\ref{control1}), we get
\begin{eqnarray}\label{controlp3}
&&C\!\! \int_{s_1}^{s_2}\!\!e^{-s}\!\!\int_{B}\!\!|w||y.\grad w|{\rho}{\mathrm{d}}y{\mathrm{d}}s\nonumber\\
&\le& C \int_{s_1}^{s_2}e^{-s}\int_{B}w^2{\rho}{\mathrm{d}}y{\mathrm{d}}s
+C e^{-s_1} \int_{s_1}^{s_2}\int_{B}|\grad w|^2(1-|y|^2){\rho}{\mathrm{d}}y{\mathrm{d}}s\nonumber\\
&\le& C e^{-s_1} \int_{s_1}^{s_2}\int_{B}w^2{\rho}{\mathrm{d}}y{\mathrm{d}}s
+K_1 e^{-s_1} +C e^{-s_1}
\int_{s_1}^{s_2}\int_{B}| w|^{p+1}\rho{\mathrm{d}}y{\mathrm{d}}s\qquad\nonumber\\
&\le& \frac{K_1}{\varepsilon}
+(K_1 \varepsilon +C e^{-s_1})
\int_{s_1}^{s_2}\int_{B}| w|^{p+1}\rho{\mathrm{d}}y{\mathrm{d}}s.\qquad
\end{eqnarray}
%We can use
Using (\ref{I12}), (\ref{CC0}), (\ref{CC1}), (\ref{CC2}) and (\ref{controlp3}), we obtain
\begin{eqnarray}\label{A1}
|A_1|&\le& \frac{K_1}{\varepsilon} +(K_1\varepsilon+
C  e^{-s_1})
\int_{s_1}^{s_2}\!\!\!\int_{B}\!| w|^{p+1}\rho{\mathrm{d}}y{\mathrm{d}}s.
\end{eqnarray}
Similarly, we deduce by (\ref{cor01}) and (\ref{control1})  that
\begin{eqnarray}\label{A2}
|A_2|
&\le&
C\int_{s_1}^{s_2} \!\!\int_{B}\!\!|ww_s|{\rho}{\mathrm{d}}y{\mathrm{d}}s \le\int_{s_1}^{s_2} \!\!\int_{B}\!\!w_s^2{\rho}{\mathrm{d}}y{\mathrm{d}}s+
C\int_{s_1}^{s_2}\!\! \int_{B}\!\!w^2{\rho}{\mathrm{d}}y{\mathrm{d}}s\nonumber\\
&\le &
K_1+
C \sup_{s\in [s_1,s_2]}\int_{B}\!\!w^2{\rho}{\mathrm{d}}y\le
\frac{K_1}{\varepsilon} +K_1\varepsilon
\int_{s_1}^{s_2}\!\!\!\int_{B}\!| w|^{p+1}\rho{\mathrm{d}}y{\mathrm{d}}s.\qquad
\end{eqnarray}
%On the other hand,
Finally, by (\ref{I100}), we obtain
\begin{eqnarray}\label{A3}
|A_3|+|A_4|&\le&  C\int_{s_1}^{s_2} \!\!e^{-\frac{2(p-q)s}{p-1}}{\mathrm{d}}s+ C\int_{s_1}^{s_2} \!\!e^{-\frac{2(p-q)s}{p-1}}\int_{B}|w|^{p+1}\rho {\mathrm{d}}y{\mathrm{d}}s\nonumber\\
&\le&  C+ Ce^{-\gamma s_1}\int_{s_1}^{s_2} \!\!\int_{B}|w|^{p+1}\rho {\mathrm{d}}y{\mathrm{d}}s.
\end{eqnarray}
This concludes the proof of Lemma 3.3 and Proposition 3.2 too.
\Box

\no Since the derivation of Theorem 1.2  from Proposition 3.2  is the
same as in the non perturbed case treated in \cite{MZ3} (up to some very minor changes),
this concludes the proof of Theorem 1.2.
\Box

\noindent{\bf Address}:\\
Universit\'e de Tunis El-Manar, Facult\'e des Sciences de Tunis, D\'epartement de math\'ematiques, Campus Universitaire 1060,
 Tunis, Tunisia.\\
\vspace{-7mm}
\begin{verbatim}
e-mail: ma.hamza@fst.rnu.tn
\end{verbatim}
Universit\'e Paris 13, Institut Galil\'ee,
Laboratoire Analyse, G\'eom\'etrie et Applications, CNRS UMR 7539,
99 avenue J.B. Cl\'ement, 93430 Villetaneuse, France.\\
\vspace{-7mm}
\begin{verbatim}
e-mail: Hatem.Zaag@univ-paris13.fr
\end{verbatim}
\end{document}